\documentclass[11pt,twoside]{article}
\usepackage{amsfonts}
\begin{document}
\centerline{\bf A general formula in Additive Number Theory}
\bigskip \centerline{Constantin M. Petridi,}
\centerline{11 Apollonos Street, 151 24 Maroussi (Athens),
Greece, cpetridi@hotmail.com} \centerline{In part - collaboration
with} \centerline{Peter B, Krikelis, Mathematics Dep., Athens
University,} \centerline{e-mail pkrikel@math.uoa.gr} \bigskip

\begin{abstract} We reduce the principal problem of Additive Number
Theory of whether an infinite sequence of integers constitutes a
finite basis for the integers to a Diophantine problem involving
the difference set of the sequence, by proving a formula
connecting the respective number of solutions.
\end{abstract}
\bigskip
\par Classical Additive Number Theory in ${\mathbb Z}$ investigates the existence of
a finite integer $\theta$ such that for a given infinite sequence
of increasing non-negative integers $a_{\nu}=f(\nu)$, $\nu =
1,2,\ldots$ the Diophantine equation $$f(x_1)+f(x_2)+\ldots
+f(x_{\theta})= n$$ has solutions $x_i \ge 1$ for any $n \ge 0$,
and determines, if possible, their number, asymptotically or
otherwise (Goldbach,  Waring, Hilbert, Hardy-Littlewood,
Vinogradov, Erd{\"o}s, \ldots ). Our approach is different.

\par
We denote by $A(n,N,\theta)= A(n.\theta)$ the number of solutions
in integers $x_i \ge 0 $ of the Diophantine system:
$$ a_1x_1+a_2x_2+\ldots+a_Nx_N = n \eqno(1)$$
$$ x_1+x_2+\ldots+x_N = \theta \hskip6pt . \eqno(1') $$

  $A(n.\theta)$
expresses in how many ways $n$ is a sum of $\theta $ integers
taken from the set $\lbrace a_1, \ldots,a_N \rbrace .$
\par
The generating function for the $A(n,\theta),\hskip6pt n = 0, 1,
\ldots ,\hskip6pt \theta = 0, 1, \ldots, \hskip6pt$ is:
$$\sum_{n=0\ldots\infty\atop \theta = 0 \ldots\infty}A(n,\theta)\hskip3pt x^ny^\theta
= \prod_{ \nu =1}^{N}{1 \over {1-x^{a_\nu}y}}\quad,$$ as seen by
expanding all the $(1-x^{a_\nu}y)^{-1}$ formally in power series
and, after multiplication, collecting terms with equal exponents
in $x$ and $y$, respectively.
\par
Considering $x$ and $y$ as complex variables we have by Cauchy's
theorem for the coefficients of power series of several complex
variables:
$$A(n,\theta)={1 \over {(2\pi}i)^2}{\int_{c_{1}}}{\int_{c_{2}}}\hskip3pt
{1 \over x^{n+1}y^{\theta +1}}\hskip3pt{\prod_{ \nu =1}^{N}}{1
\over {1-x^{a_\nu}y}}\hskip3ptdxdy\hskip3pt , \quad\eqno(2) $$
where the integrals are taken over the circles $|x|=c_1 <
1$,\hskip3pt $|y|=c_2 < 1$, respectively.
\par
The partial fraction expansion with regard to $y$ of $\hskip3pt
{\prod_{ \nu =1}^{N}}({1-x^{a_\nu}y})^{-1}$ gives
$${\prod_{ \nu =1}^{N}}{1
\over {1-x^{a_\nu}y}} = \sum_{\nu=1}^N{1 \over
y^{N^{-1}}L'(x^{a_\nu})(1-x^{a_\nu}y)}\hskip6pt ,$$ where $L'(t)$
is the derivative of $L(t)=\prod_{\nu=1}^N(t-x^{a_\nu})$.
\par
Inserting in (2) and expanding ${(1-x^{a_\nu}y})^{-1}$ in power
series of $y$ within the circle $c_2$ we get successively:
$$A(n,\theta)={1\over (2\pi i)^2}{\int_{c_{1}}}{\int_{c_{2}}}\hskip3pt {1\over x^{n+1}y^{\theta
+1}}\hskip3pt {\biggl\lbrace \sum_{\nu=1}^N {1 \over
y^{N^{-1}}L'(x^{a_\nu})(1-x^{a_\nu}y)}  } \biggr\rbrace dxdy
\quad$$
$$\qquad = {1\over 2\pi i}{\int_{c_{1}}}\hskip3pt {1 \over x^{n+1}}\hskip3pt \biggl\lbrace {1\over 2\pi i}{\int_{c_{2}}}\hskip3pt
{1 \over y^{\theta +N}}\hskip3pt {\sum_{\nu=1}^N }{1 \over
L'(x^{a_\nu}) \bigl( 1-x^{a_{\nu}}y\bigr )}\hskip3ptdy
\biggr\rbrace dx$$
$$\qquad\qquad \qquad ={1\over 2\pi i}{\int_{c_{1}}}\hskip3pt {1 \over x^{n+1}} \hskip3pt \biggl\lbrace {\sum_{\nu=1}^N }
\hskip3pt {1\over 2\pi i}{\int_{c_{2}}}\hskip3pt {1 \over
y^{\theta +N}L'(x^{a_{\nu}})} \bigl( 1+x^{a_{\nu}}y+ \ldots
\hskip6pt \bigr )dy \biggr\rbrace dx \quad . $$ Integrating over
$|y|=c_2$ we obtain by Cauchy's theorem:
$$A(n,\theta)={1\over 2\pi i}{\int_{c_{1}}}\hskip3pt {1 \over x^{n+1}} \biggl\lbrace
{\sum_{\nu=1}^N }\hskip3pt {x^{a_{\nu}(\theta +N-1)} \over
L'(x^{a_{\nu}})} \biggr\rbrace dx  $$
$$ \qquad\qquad \qquad ={\sum_{\nu=1}^N }\hskip3pt{1\over 2\pi i}{\int_{c_{1}}}\hskip3pt{1 \over x^{n+1-a_{\nu}(\theta
+N-1)}}{1 \over L'(x^{a_{\nu}})}\hskip3pt dx .\eqno(3)$$ In order
to expand $1/L'(x^{a_{\nu}}) $ in a power series of $x$ we
transform $L'(x^{a_{\nu}})$ as follows:
$$L'(x^{a_{\nu}})=(x^{a_{\nu}}-x^a_1)\ldots(x^{a_{\nu}}-x^{a_{\nu
-1}})(x^{a_{\nu}}-x^{a_{\nu +1}})\ldots(x^{a_{\nu}}-x^{a_N} )$$
$$ = x^{a_1+\ldots+a_{\nu
-1}}(x^{a_{\nu}-a_1}-1)\ldots(x^{a_{\nu}-a_{\nu-1}}-1)
x^{a_{\nu}(N-\nu)}(1-x^{a_{\nu
+1}-a_{\nu}})\ldots(1-x^{a_{N}-a_{\nu}}) $$
$$ =(-1)^{\nu -1}x^{a_1+\ldots+a_{\nu
-1}+a_{\nu}(N-\nu)}(1-x^{a_{\nu}-a_1})\ldots(1-x^{a_{\nu}-a_{\nu
-1}})(1-x^{a_{\nu +1}-a_{\nu}})\ldots(1-x^{a_N-a_{\nu}}). $$
Inserting in (3) for $L'(x^{a_{\nu}})$ their above expressions and
after calculations in the exponents, we obtain:
$$ {A(n,\theta)={\sum_{\nu=1}^N }\hskip3pt {(-1)^{\nu -1} \over {2\pi i}}{\int_{c_{1}}}
{1 \over x^{n+1+a_1+\ldots+a_{\nu}-(\theta
+\nu)a_{\nu}}}P_{\nu}(x) dx } \eqno(4)
$$ where
$$P_{\nu}(x)={1 \over
(1-x^{a_{\nu}-a_1})\ldots(1-x^{a_{\nu}-a_{\nu -1}})(1-x^{a_{\nu
+1}-a_{\nu}})\ldots(1-x^{a_N-a_{\nu}})} \hskip6pt .$$ Since all
exponents $a_{\nu}-a_{\mu}$ are positive, we can expand for
$|x|\le c_1$ the factors $P_{\nu}(x)$ in power series of $x$:
$$ P_{\nu}(x) = {\sum_{\lambda =0}^{\infty}} B_{\nu} (\lambda)x^{\lambda }
$$ where $B_{\nu}(\lambda )$ are respectively the number of non-negative integer solutions
 of the linear Diophantine equations
$$(a_{\nu}-a_1)x_1+ \ldots +(a_{\nu}-a_{\nu -1})x_{\nu -1}+
(a_{\nu +1}-a_{\nu})x_{\nu}+ \ldots +(a_{N}-a_{N-1})x_{N-1}=
\lambda \hskip6pt .$$

\par
Substituting in (4) the $P_{\nu}(x)$ by their respective power
series and using again Caushy's theorem we finally arrive at
$$ A(n, \theta)=\sum_{\nu =1}^N \hskip3pt (-1)^{\nu-1}B_{\nu}(s_{\nu} )\hskip6pt,\eqno(5)$$
with $s_{\nu}=n+ {\sum_{i=1}^{\nu}}a_i-(\nu +\theta)a_{\nu}$ and
$B_{\nu}(s_{\nu})$ respectively, the number of solutions of each
of the following linear Diophantine equations
$$(a_{\nu}-a_1)x_1+ \ldots +(a_{\nu}-a_{\nu -1})x_{\nu -1}+
(a_{\nu +1}-a_{\nu})x_{\nu}+ \ldots +(a_{N}-a_{\nu -1})x_{N-1}=
s_{\nu} \eqno(6)$$
$$\nu = 1, \ldots, N  \hskip6pt .$$
\par
This formula reduces the investigation of the number of solutions
of the initial system to that of the number of solutions of $N$
linear Diophantine equations involving the difference sets
(positive) $\lbrace a_{\nu}-a_{\mu} \rbrace $ and the numbers
$s_{\nu}$.
\par
Geometrically speaking this means that the number of {\it
Gitterpunkte} in the intersection of the two $N$-dimensional
planes (1) and (1') in the positive quadrant $x_i \ge 0$ is equal
to the alternate sum of the number of {\it Gitterpunkte} of the
$(N-1)$-dimensional planes (6) in the same quadrant.
\par
As standard examples, theorems and conjectures we may cite
$a_{\nu}=(\nu -1)^2$, $\theta = 4$ (Lagrange), $a_{\nu}=(\nu
-1)^k$ (Waring), $a_{\nu}=\nu $-th prime, $\theta = 2$, $n$ even
(Goldbach), $a_{\nu}=\nu ^p$, $\theta = 2$, $n=n_1^p$, $p \ge 3$
(Fermat).
\par
Obviously in order to attack the problem for a given sequence we
have to take into account that $N$ is linked to $n$ by a function
$N(n)$ (ex.g. for Lagrange $N(n)=\lbrack n^{1\over 2} \rbrack$).
This complicates the matter but still an ad hoc suggestion would
be to approximate the $B_{\nu}(s_{\nu})$ as follows
 $$B_{\nu}(s_{\nu})\sim {s_{\nu}^{N-2} \over
(N-2)!(a_{\nu}-a_1) \cdots (a_N-a_{\nu})} \hskip6pt , $$ which is
valid for $n \rightarrow \infty$ but fixed $N$ (Polya-Szeg{\"o},
Aufgaben und Lehrs{\"a}tze aus der Analysis I; loosing no
generality the differences $a_{\nu}-a_{\mu}$ can be assumed free
of common divisors $>1$).
\par
Summing over $\nu$ (we write $N$ for short of $N(n)$) and
reverting again to the polynomials $L(t)$, written now as
$L_N(t)$, we obtain from (5)
$${1 \over (N-2)!}\sum_{\nu=1}^{N}\hskip3pt {s_{\nu}^{N-2} \over L_N'(a_{\nu})} \hskip6pt ,\eqno(7)$$
as a plausible heuristic estimate of $A(n, \theta)$ for $n
\rightarrow \infty$.
\par
The behaviour of the expressions involving $n$ and $\hskip4pt
\theta$:

\begin{center}
\begin{tabular}{ll}
 $\displaystyle{{n+ {\sum_{i=1}^{\nu}}a_i-(\nu +\theta)a_{\nu}} \over {|a_{\nu}-a_{\mu}|}}$ & The cutting
points of the planes (6) with the coordinate axes,\\
 & \\
 $\displaystyle{{ n+ {\sum_{i=1}^{\nu}}a_i-(\nu +\theta)a_{\nu}} \over
{\sqrt{\sum_{\nu=1}^N(a_{\nu}-a_{\mu})^2}}} $& The distances of
the planes (6) from the origin,
\end{tabular}
\end{center}

would play, we believe, a decisive role in any such attempt.
\par
As to its form the sum (7) bears a striking resemblance to the
sums
$$\sum_{\nu =1}^N{a_{\nu}^t \over {L'_N(a_{\nu})}}$$
encountered in Lagrange interpolation with $a_{\nu}$ replaced by
$s_{\nu}$ in the numerator. As known these expressions considered
as functions of the exponent $t$ are equal to:
$$   \left\{\begin{tabular}{l@{\quad}llcl}
                          $\displaystyle {(-1)^{N-1} \over {a_1 \cdots a_N}}
 \displaystyle \sum_{\sum i=-t-1}\hskip3pt {1 \over {a_1^{i_1} \cdots
 a_N^{i_N}}}$ & & for & \quad\hskip3pt $t \le -1$  \\
                     \qquad\quad\quad $0$ & & for & \quad $0 \le t \le N-2$\\
\qquad\quad\quad $1$ & & for & \quad \quad\quad $t=N-1$ \\
$\displaystyle \sum_{\sum i=t-N+1} \hskip3pt a_1^{i_1} \cdots
a_N^{i_N}$ & &for & $N \le t $\qquad\quad\quad
                          \end{tabular} \right.    $$

Above facts may prove, eventually, useful in further developments.
\end{document}